\documentclass[draft]{amsart}

\usepackage{amssymb,latexsym}
\def\C{\mathbb C}
\def\R{\mathbb R}
\def\N{\mathbb N}
\def\Z{\mathbb Z}

\newtheorem{thm}{Theorem}[section]
\newtheorem{lem}{Lemma}[section]

\begin{document}
%\sffamily
\title{Non-real zeros of derivatives of real meromorphic functions}
\author{J.K. Langley}
\address
{School of Mathematical Sciences
\\University of Nottingham
\\Nottingham NG7 2RD}
\email{jkl@maths.nottingham.ac.uk}
\thanks{Research supported by Engineering and Physical Sciences Research Council grant EP/D065321/1}

\maketitle
\begin{abstract}
The main result of the paper determines all real meromorphic functions $f$ of finite order in the plane such that $f'$ has finitely many zeros while $f$ and $f^{(k)}$, for some $k \geq 2$, have finitely many non-real zeros.
\end{abstract}

\section{Introduction}

This paper concerns non-real zeros of the derivatives of real meromorphic functions in the plane: here a 
meromorphic function is called real if it maps 
$\R$ into $\R \cup \{ \infty \}$. In the setting of real entire functions \cite{HelW1,HelW2,SS} the class $V_{2p}$ is defined for $p \geq 0$ to consist of all entire functions
$f(z) = g(z) \exp( -a z^{2p+2} ) $,
where $a \geq 0$ is real and $g$ is a real entire function 
with real zeros of 
genus at most $2p+1$ \cite[p.29]{Hay2}. It is well known \cite{Lag}
that $V_0 $ coincides with the Laguerre-P\'olya class $LP$ of
entire functions which are locally uniform limits of real polynomials
with real zeros. With the notation $V_{-2} = \emptyset$ the class  $U_{2p}^*$ may then be defined for $p \geq 0$ as the set of entire functions $f = Ph$, where 
$h \in  V_{2p} \setminus V_{2p-2}$ and $P$ is a real polynomial without
real zeros \cite{EdwH}; thus each real
entire function of finite order with finitely many non-real zeros
belongs to $U_{2p}^*$ for some $p \geq 0$.
The following results, in which all counts of zeros are with respect to multiplicity, established conjectures of Wiman \cite{Al1,Al2}
and P\'olya \cite{Polya43} respectively.

\begin{thm}[\cite{EdwH,SS}]\label{thmC}
Let $p \in \N$ and let $f \in U_{2p}^*$. Then $f''$
has at least $2p$ non-real zeros.
\end{thm}

\begin{thm}[\cite{BEpolya}]\label{thmB}
Let $p$ be a positive integer and let $f \in U_{2p}^*$. Then the number
of non-real zeros of the $k$th derivative
$f^{(k)}$ tends to infinity with $k$.
\end{thm}
\begin{thm}[\cite{BEL,lajda}]\label{thmA}
If $f$ is
a real entire function of infinite order 
then $ff^{(k)}$ has infinitely many
non-real zeros, for every $k \geq 2$.
\end{thm}

For real meromorphic functions with poles rather less is known. All meromorphic functions $f$ in the plane for which all derivatives $f^{(k)} $ ($k \geq 0$) have only real zeros were determined by Hinkkanen in a series of papers \cite{Hin1,Hin2,Hin3}: such functions have at most two distinct poles, by the P\'olya shire theorem \cite[Theorem 3.6, p.63]{Hay2}. Functions with real poles, for which some of the derivatives have only real zeros, were treated in a number of papers including \cite{HelW3,HSW,rossireal}. In particular the following theorem was proved in \cite{HSW}.

\begin{thm}\label{thmHSW}
 Let $f$ be a real meromorphic function in the plane with only real zeros and poles (and at least one of
 each). Assume that $f'$ has no zeros and that $f''$ has only real zeros. Then $f$ has one of the forms
$$
A \tan (az+b) + B, \quad \frac{az+b}{cz+d} , \quad A \cdot \frac{(az+b)^2-1}{(az+b)^2} ,
$$
where $A, B, a, b, c, d$ are real numbers.
\end{thm}

The following result will be proved: here the reality of all but finitely many poles of $f$ is no longer an assumption but turns out to be a conclusion, and the second derivative is shown to be exceptional. 

\begin{thm}\label{thm1}
Let $f$ be a real meromorphic function in the plane, not of the form 
$f = Re^P$ with $R$ a rational function and $P$ a polynomial. Let $k \geq 2$ be an integer.
Assume that:\\
(i) all but finitely many zeros of $f$ and $f^{(k)}$ are real;\\
(ii) the first derivative $f'$ has finitely many zeros;\\
(iii) 
%either $k = 2$ and $f$ has finite order, or $k \geq 3$ and 
there exists $M \in (0, \infty)$
such that if $\zeta $ is a pole of $f$ of multiplicity $m_\zeta$ then
\begin{equation}
m_\zeta \leq M + |\zeta|^M ;
\label{multest}
\end{equation}
(iv) if $k=2$ then $f'/f$ has finite order.

Then $f$ satisfies
%there exist $c$, $A$ and $R$ such that
%%positive real number $c$, a complex number $A$ and a rational function $R$
%with $|R(x)| = 1$ for all $x \in \R$,  and $A \in \C \setminus \R$ 
\begin{equation}
 f(z) = \frac{ R(z)e^{icz} -1}{A R(z) e^{icz} - \overline{A}} \, , \quad
\hbox{where $c \in (0, \infty) $, $ A \in \C \setminus \R$,} 
\label{conc1}
\end{equation}
and
\begin{equation}
\hbox{$R$ is a rational function with 
$|R(x)| = 1$ for all $x \in \R$.} 
\label{conc}
\end{equation}
Moreover, $k=2$ and all but finitely many poles of $f$ are real.

Conversely, if $f$ is given by (\ref{conc1}) and (\ref{conc}) then $f$ satisfies (i) and (ii) with $k=2$.
\end{thm}
If the function $f$ is given by 
$f = Re^P$ with $R$ a rational function and $P$ a polynomial then obviously $f$ and all its derivatives have finitely many zeros.
%Lemma \ref{lemelem2} will show that if the function $f$ satisfies the conclusions of Theorem \ref{thm1} %then it satisfies the hypotheses with $k=2$.
It is not clear whether the assumptions (iii) and (iv) are really necessary in Theorem \ref{thm1}, but they are required for the proof presented below. 
In the case $k\geq3$ it will be proved in \S\ref{fgrowth} that if a real meromorphic function $f$ satisfies (i) and (ii) 
%and (\ref{multest}) 
then $f'/f$ has
finite order: this will use a method of Frank \cite{Fra1}, which is not available for $k=2$. On the other hand, if $f$ itself has finite order then (\ref{multest}) holds since, with the standard notation of Nevanlinna theory \cite{Hay2}, 
$$
n(r, f) = O( N(2r, f) ) = O( T(2r, f) ).
$$
%The simple example $f(z) = \tan z$ gives
%$$
%f''(z) = 2 \sec^2 z \tan z, \quad
%f'''(z) = 6 \sec^4 z - 4 \sec^2 z , \quad f^{(iv)}(z) = 
%(24 \sec^4 z - 8 \sec^2 z) \tan z .
%$$
%Thus $f$ satisfies the hypotheses of Theorem \ref{thm1} with $k=2$, but not $k=3$. Indeed $f'''$ has no %real zeros (cf. \cite[Theorem 2]{HelW3}), but a further ......

\section{Preliminaries}
%\begin{definition}\label{def1}
For $a \in \C$ and $r > 0$ set
$D(a, r) = \{ z \in \C : | z- a | < r \}$ and correspondingly let $S(a, r) = \{ z \in \C : | z- a | = r \}$.
%and 
%$$
%A(r, \infty ) = \{ z \in \C \cup \{ \infty \}  : r < |z| \leq \infty \} .
%$$
%Further, set
%\begin{equation}
%H^+ = \{ z \in \C : {\rm Im} \, z > 0 \}, \quad
%D^+(0, r) = D(0, r) \cap H^+, \quad
%A^+(r, \infty ) = A(r, \infty ) \cap H^+ .
%\label{A1}
%\end{equation}
%\end{definition}
The following lemma is standard \cite[pp.116-7]{Tsuji}.
\begin{lem}[\cite{Tsuji}]\label{lemhm}
Let $u$ be a non-constant continuous subharmonic function in the plane.
For $r > 0$
let $\theta^*(r)$ be the angular measure of that subset of $S(0, r)$
on which $u(z) > 0$, except that $\theta^*(r) = \infty $ if 
$u(z) > 0$ on the whole circle $S(0, r)$. Then, if $r \leq R/4$ and $r$ is sufficiently large,
\begin{equation}
B(r, u) = \max \{ u(z) : |z| = r \}
%\leq \frac3{2 \pi} 
%\int_0^{2 \pi} \max \{ u( 2r e^{it} ) , 0 \} dt 
%\label{p2}
%\end{equation}
%and,  
%\begin{equation}
%B(r, u) 
\leq 9 \sqrt{2} B(R, u) \exp \left( - \pi \int_{2r}^{R/2} \frac{ds}{s \theta^*(s)} \right).
\label{p3}
\end{equation}
\end{lem}
%The inequality
%\ref{p2}) follows from Poisson's formula, and (\ref{p3}) from
%a standard application of a well known estimate for harmonic measure
%\cite[pp.116-7]{Tsuji}.

Next, let the function $g$ be meromorphic
in a domain containing the closed upper half-plane $\overline{H} =
\{ z \in \C : \mathrm{Im} \, z \geq  0 \}$.
For $t \geq 1$ let
$\mathfrak{n} (t, g) $ be the number of poles of $g$, counting
multiplicity, in $\{ z \in \C:|z-it/2|\leq t/2, |z| \geq 1\}$.
The Tsuji characteristic $\mathfrak{T} (r, g)$ \cite{GO,LeO,Tsuji0} is defined for $r \geq 1$ by 
$\mathfrak{T} (r, g)  =
\mathfrak{m} (r, g)  +
\mathfrak{N} (r, g)$, where 
$$
\mathfrak{m} (r, g) =
\frac1{2 \pi} \int_{ \sin^{-1} (1/r)}^{ \pi - \sin^{-1} (1/r)}
\frac{ \log^+ | g(r \sin \theta e^{i \theta } )|}
{r \sin^2 \theta } d \theta
\quad \hbox{and} \quad \mathfrak{N} (r, g) = \int_1^r
\frac{ \mathfrak{n} (t, g) }{t^2} dt.
$$

\begin{lem}[\cite{LeO}]\label{lemLO}
Let the function $g$ be meromorphic in 
%the closed upper half-plane 
$\overline{H}$ and assume that
%\begin{equation}
$\mathfrak{m}(r, g) = O( \log r )$ as $ r \to \infty $.
%\label{p4}
%\end{equation}
%where $\mathfrak{m}(r, g) $ is given by $(\ref{tsujim})$.
Then, as $R \to \infty$, 
$$
\int_R^\infty  \int_0^\pi \frac{ \log^+ |g(r e^{i \theta  } )| }{r^3}\, d \theta \, dr 
= O\left( \frac{ \log R }{R} \right) .
$$
\end{lem}

The following theorem was proved in \cite{BLa,Schwick}
using the rescaling method \cite{Zalc}. 
\begin{thm}[\cite{BLa,Schwick}]\label{norfam}
Let $k \geq 2$ and let ${\mathcal F}$ be a family of functions meromorphic
on a plane domain $D$ such that $ff^{(k)}$ has no zeros in $D$, for each
$f \in {\mathcal F}$. Then the family
$\{ f'/f : f \in {\mathcal F} \}$ is normal on $D$.
\end{thm}

The proof of Theorem \ref{thm1} will require well known results of Edrei-Fuchs \cite{EF2} and Hayman \cite{Hay3} respectively. 

\begin{lem}[\cite{EF2}]\label{lemEF}
Let $1 < r < R < \infty $ and let the function $g$ be meromorphic in
$|z| \leq R$. Let $I(r)$ be a subset of $[0, 2 \pi ]$ of Lebesgue
measure $\mu (r)$. Then
$$
\frac1{2 \pi} \int_{I(r)} 
\log^+ | g(re^{i \theta })| d \theta \leq 
\frac{ 11 R \mu (r)}{R - r} \left( 1 + \log^+ \frac1{\mu(r)} \right)
T(R, g).
$$
\end{lem}

\begin{lem}[\cite{Hay3}]\label{Haylem}
Let 
$S(r)$ be an unbounded positive non-decreasing
function on $[r_0, \infty )$, continuous from the right, of
order $\rho$. Let $A > 1, B > 1$. Then
$$
\overline{{\rm logdens }} \, G \leq \rho \left( \frac{\log A}{\log B} \right), 
\quad \hbox{where} \quad 
G = \{ r \geq r_0 : S(Ar) \geq B S(r) \}. 
$$
\end{lem}

The next lemma \cite{EL,ripstall} has found widespread applications in function theory and complex dynamics.

\begin{lem}[\cite{EL,ripstall}]\label{BELlem}
Let the function $g$ be transcendental and meromorphic in the plane such that the set
of finite singular values of the inverse function $g^{-1}$ is bounded.
Then there exist $L > 0$ and $M > 0$ such that
$$
\left|  \frac{z_0 g'(z_0) }{g (z_0) } \right|
\geq
C \log^+ \left|\frac{g(z_0)}{M}\right|
\quad \hbox{for} \quad |z_0| > L,
$$
where $C$ is a positive absolute constant, in particular
independent of $g, L$ and $M$.
\end{lem}

\section{An elementary lemma and some consequences}\label{elem}

\begin{lem}
 \label{lemelem}
Let 
$$
g(z) = \frac1{1 - e^z} .
$$
%Then all zeros of $g''$ lie on the imaginary axis, and 
Then for $k \geq 3$ the function $g^{(k)}$ has infinitely many zeros off the imaginary axis.
\end{lem}
\begin{proof} With the notation $X = e^z$, 
%elementary computation gives
%$$g'(z) = \frac{X}{(1-X)^2} , \quad g''(z) = \frac{X^2 + X}{(1-X)^3} .$$
%and the first assertion of the lemma is obvious. 
it will be proved by induction that
\begin{equation}
 \label{gk}
g^{(k)}(z) = \frac{P_k(X) }{(1-X)^{k+1}} =\frac{X^k + A_{k-1} X^{k-1} + \ldots }{(1-X)^{k+1}},
\end{equation}
in which the numerator $P_k(X)$ is a monic polynomial in $X$ of degree $k$ with constant coefficients,  and the coefficient of $X^{k-1}$ is denoted by $A_{k-1}$. Evidently (\ref{gk}) is true for $k=1$ and $k=2$, with $A_0 = 0$ and $A_1 = 1$. Assuming that $k \geq 2$ and that (\ref{gk}) holds, differentiation gives
%$$g^{(k+1)}(z) = \frac{P_{k+1}(X) }{(1-X)^{k+2}}$$in which
\begin{eqnarray*}
P_{k+1}(X) &=&
(1-X)^{k+2} g^{(k+1)}(z)\\
&=& (kX^k + (k-1)A_{k-1}X^{k-1} + \ldots )(1-X) + \\
& & + (k+1)X(X^k + A_{k-1}X^{k-1} + \ldots )\\
&=& X^{k+1} + X^k ( k - (k-1)A_{k-1} + (k+1) A_{k-1}) + \ldots \\
&=&
X^{k+1} + (k + 2 A_{k-1} ) X^k + \ldots .
\end{eqnarray*}
This proves (\ref{gk}) with $k$ replaced by $k+1$, and gives in addition the recurrence relation 
$A_k = k + 2 A_{k-1}$. Since $A_1 = 1$ it follows at once that $A_k \geq k+2$ for $k \geq 2$. 
But then for $k \geq 3$ the sum of the roots of $P_k(X)$ is $-A_{k-1} \leq -k-1$, and so these roots cannot all have modulus $1$. 
\end{proof}

\begin{lem}\label{lemelem2}
 Let the function $f$ be given by (\ref{conc1}) and (\ref{conc}).
%, in which $c$ is a positive real number, $A$ is a 
%non-real complex number, and $R$ is a rational function with $|R(x)| = 1$ for all $x \in \R$. 
Then all but finitely many zeros of $f''$ and poles of $f$ are real, and for $k \geq 3$ the function $f^{(k)}$ has infinitely many non-real zeros.
\end{lem}
\begin{proof}
The fact that all but finitely many poles of $f$ and zeros of $f''$ are real is proved in \cite{nicks}, but a slightly different argument is given here for completeness. First, it is evident from (\ref{conc1}) and (\ref{conc}) that there exist constants $D, E$ and
rational functions $S,  U, V$ with
\begin{equation}
 \label{elem1}
f(z) = D + \frac{E}{1 - S(z)e^{icz} } 
\quad \hbox{and} \quad 
f'' (z) = \frac{U(z)e^{icz} + V(z)e^{2icz}}{(1 - S(z)e^{icz})^3} ,
\end{equation}
where $|S(x)| = 1$ for $x \in \R$. 
If $z_k$ is a pole of $f$ with $|z_k|$ large then 
$e^{ic z_k} = C + o(1)$, where $C \cdot S(\infty) = 1$ and so $|C| = 1$. 
Hence $z_k$ lies near a zero $x_j$ of $e^{icz} - C$, and $x_j$ is real. Moreover if $|x_j|$ is large enough then Rouch\'e's theorem gives precisely one pole of $f$ near $x_j$. Since $f$ is real, so is $z_k$. 

Next, let $x \in \R$. Since $f''(x)$ is real and $|S(x)| = 1$, the representation (\ref{elem1}) gives
$$
f''(x) = \overline{f''(x)} = \frac{\overline{U(x)}e^{-icx} + \overline{V(x)}e^{-2icx} }
{(1- \overline{S(x)}e^{-icx})^3} = 
- \frac{S(x)^3\overline{U(x)}e^{2icx} + S(x)^3 \overline{V(x)}e^{icx} }{(1 - S(x)e^{icx})^3}
$$
so that $S(x)^3\overline{U(x)} = - V(x)$ and $|V(x)/U(x)| = 1$. The same argument as for the poles now shows that all but finitely many zeros of $f''$ are real (this may also be proved using the Levin-Ostrovskii representation \cite{LeO} for $-f''/f'$, since the real function $1/f'$ has finitely many poles and finitely many non-real zeros). 

To prove the last assertion let $k \geq 3$. Write $S(\infty) = e^{id} $ for some $d \in \R$ and
$$
 h(z) = \frac1{1 - S(z)e^{icz} } = \frac1{1 - T(z)e^{i(cz+d)}}, \quad
H(z) = \frac1{1 - e^{i(cz+d)}} = g(i(cz+d)) ,
$$
in which the function $T$ is rational with $T(\infty) = 1$, and $g$ is as in Lemma \ref{lemelem}.
By (\ref{elem1}) it suffices to prove that $h^{(k)}$ has infinitely many non-real zeros. Lemma \ref{lemelem} gives a non-real zero $w$ of $H^{(k)}$ and there exist small positive $s, t$ such that
$$
|H^{(k)}(z)| \geq s \quad \hbox{and} \quad s \leq |1 - e^{i(cz+d)} | \leq \frac1s 
\quad \hbox{for $t \leq |z - w| \leq 3 t$.}
$$
Let $n$ be a large positive integer. Then the periodicity of $H$ yields
$$h(z) = \frac1{1 - e^{i(cz+d)}(1+o(1) )} =\frac{1 + o(1)}{1 - e^{i(cz+d)}} = H(z) + o(1)
$$
for $t \leq |z - (w+n 2\pi/c)| \leq 3 t$. This implies that
$$h^{(k)} (z) = H^{(k)}(z) + o(1) = H^{(k)}(z) (1 + o(1) ) $$
on the circle  $S(w + n 2 \pi/c, 2t)$, and so $h^{(k)}$ has a zero in $D(w + n 2 \pi /c, 2t)$ by Rouch\'e's theorem.
\end{proof}

\section{Proof of Theorem \ref{thm1}: the first part}\label{fgrowth}
Let the integer $k $ and the function $f$ be as in the statement of Theorem \ref{thm1}. Since
$f$ is not of the form $f = Re^P$ with $R$ a rational function and $P$ a polynomial, the logarithmic derivative $L =  f'/f$ is transcendental. 
The first task is to show that $f'/f$ has finite order.
For $k=2$ this is true by hypothesis (iv), and for the case $k\geq3$ the fact that $f$ and $f^{(k)}$ have finitely many non-real zeros leads at once to the following lemma, which is
proved using a method of Frank \cite{Fra1,FHP} (see also \cite{BLa,FH,FL2}), 
but with the Nevanlinna characteristic replaced by that of Tsuji. 
%This method is, however, not available for $k = 2$. 

\begin{lem}[\cite{Fra1}]\label{lemA1}
Assume that $k \geq 3$. Then the Tsuji characteristic of
$L = f'/f$ satisfies
\begin{equation}
%\mathfrak{m}(r, 1/L) \leq
\mathfrak{T}(r, L) = O( \log r) \quad\hbox{as}\quad r \to \infty.
\label{A5}
\end{equation}
\end{lem}

\begin{lem}
 \label{lemfo1}
%Assume again that $k \geq 3$. 
The function $L = f'/f$ has finite order.
\end{lem}
\begin{proof} 
If $k=2$ there is nothing to prove, so assume that $k \geq 3$. 
It follows at once from (\ref{A5}) that
$$
\mathfrak{m}(r, 1/L) \leq
\mathfrak{T}(r, L) + O(1) = O( \log r) \quad\hbox{as}\quad r \to \infty.
$$
Hence Lemma \ref{lemLO} and the fact that $f$ is real give
$$
\int_R^\infty \frac{m (r, 1/L) }{r^3} dr 
=  O\left( \frac{ \log R }{R} \right) \quad\hbox{as} \quad R \to \infty .
$$
But $1/L = f/f'$ has finitely many poles and it now follows that, as $R \to \infty$,  
$$
\frac{T(R, 1/L)}{R^2} \leq 2 \int_R^\infty \frac{T (r, 1/L) }{r^3} dr 
\leq 2 \int_R^\infty \frac{m (r, 1/L) }{r^3} dr +  O\left( \frac{ \log R }{R^2} \right)
=  O\left( \frac{ \log R }{R} \right).
$$
\end{proof}

\begin{lem}
 \label{lemfo4}
There exists a set $E_1 \subseteq [1, \infty)$, 
of upper logarithmic density at most $1/2$, with the following property. 
To each real number $\sigma \in (0, \pi /2)$ corresponds a positive real number $N_1$ such that, for all large $r \not \in E_1$,
\begin{equation}
 \label{fo4}
\log |r L(re^{i \theta } )| \leq N_1  \quad \hbox{for $\sigma \leq \pm \theta \leq \pi - \sigma$.}
\end{equation}
\end{lem}
\begin{proof} 
Set
$$
g(z) = \frac{f(z)}{zf'(z)} .
$$
Then $g$ is transcendental of finite order, by Lemma \ref{lemfo1}. Lemma \ref{Haylem} gives positive real numbers
$B_1$ and $r_1$, with $r_1$ large, such that 
\begin{equation}
 \label{fo3}
T(2r, g) \leq B_1 T(r, g) 
\end{equation}
for all $r \geq r_1$ outside a set $E_1$ of upper logarithmic density at most $1/2$. 
For $r \geq r_1$ let $I_1(r)$ be that
subset of $[0, 2 \pi]$ on which $|g(re^{i \theta } )| \geq 1$ and let $\mu_1(r)$ denote the Lebesgue measure of $I_1(r)$. By (\ref{fo3}), Lemma \ref{lemEF} and the fact that $g$
has finitely many poles, 
there exists a positive real number $\tau$ such that $\mu_1(r) > 5 \tau$ for all $r \geq r_1$ with $r \not \in E_1$. 

Let $\sigma \in (0, \pi /2)$.
% and assume without loss of generality that $ \sigma \leq \tau$. 
Theorem \ref{norfam} implies that the family 
$$
\mathcal{G} = \{ G(z) = r f'(rz)/f(rz) : r \geq r_1, \, r  \not \in E_1 \}
$$
is normal on the domain $\{ z \in \C : 1/2 < |z| < 2 , \, 0 < \arg z < \pi \}$. By the definitions
of $g$ and $I_1(r)$ and the fact that $f$ is real there exists, for each $G \in \mathcal{G}$,
a real number $\theta_1 \in ( \tau, \pi - \tau) $ 
such that $|G(e^{i \theta_1} )| \leq 1$. A standard normal families argument then gives
$|G(e^{i \theta} )| \leq M_\sigma$ for $\sigma \leq \theta \leq \pi - \sigma$, where 
$M_\sigma$ depends only on $\sigma$, and (\ref{fo4}) follows, using again the fact that $f$ is real.
\end{proof}

\begin{lem}
 \label{lemfo2}
The function $f$ has finite order.
\end{lem}
\begin{proof} 
%Assume that $k \geq 3$, since otherwise there is nothing to prove.
Combining Lemma \ref{lemfo1} with (\ref{multest}) and the hypothesis that $f'$ has finitely many zeros 
shows that the zeros and poles of $f$ have finite exponent of convergence. Hence there exist
a real meromorphic function $\Pi$ of finite order and a real entire function $h$ such that
\begin{equation}
 \label{fo1}
f = \Pi e^h, \quad L = \frac{f'}{f} = \frac{\Pi'}{\Pi} + h' .
\end{equation}
Since $L$ and $\Pi$ have finite order, so has $h$. Assume that $h$ is transcendental: if this is not the case there is nothing to prove. Standard estimates for logarithmic derivatives \cite{Gun2} show that there
exists a positive real number $N_2 $ such that
\begin{equation}
 \label{fo2}
\log \left| \frac{\Pi'(z)}{\Pi(z)} \right| \leq N_2 \log |z| 
\end{equation}
provided that $|z|$ is large and lies outside a set $E_2 \subseteq [1, \infty)$ of finite logarithmic measure. 
Let $\sigma$ be a small positive real number and let $N_1$ and the set $E_1$ be as in Lemma \ref{lemfo4}.
For large $r$ let 
$$I_2(r)= \{ \theta \in [0, 2 \pi]: \, \log |h'(re^{i \theta } )| \geq ( N_2 + 1) \log r \}.$$
Since $h'$ is real, it follows from (\ref{fo4}), (\ref{fo1}) and 
(\ref{fo2}) that, for all large $r \not \in E_3 = E_1 \cup E_2$,
the set $I_2(r) $ may be enclosed in a set $ I_3(r) \subseteq [0, 2 \pi ]$ of
Lebesgue measure $ 4 \sigma $.
% , which has upper logarithmic density at most $1/2$. 
Applying Lemma \ref{lemEF} to the transcendental entire function $h'$ yields
\begin{eqnarray*}
T(r, h') &\leq& \frac1{ 2 \pi} \int_{I_3(r)} 
\log^+ | h'(re^{i \theta })| d \theta + O( \log r ) \\
&\leq& 176 \sigma 
\left( 1 + \log^+ \frac1{4  \sigma } \right)
T(2r, h') = \delta T(2r, h') 
\end{eqnarray*}
for all large $r \not \in E_3 = E_1 \cup E_2$ and so for all $r$ in a set
$E_4 \subseteq [1, \infty)$ of lower logarithmic density at least $1/2$. Lemma \ref{Haylem} now shows
that the order $\rho (h')$ of $h'$ satisfies 
$$
\frac12 \leq \underline{{\rm logdens }} \,E_4 \leq 
\overline{{\rm logdens }} \, E_4 \leq \rho (h') \left( \frac{\log 2}{\log 1/\delta } \right) .
$$
But $\delta$ may be made arbitrarily small by choosing $\sigma$ small enough, and this is a contradiction.
\end{proof}

\section{Asymptotic values of $f$}\label{fav}

The next step in the proof of Theorem \ref{thm1} is to show that $f$ has precisely two asymptotic values, both of them finite and non-real. Since $f$ has finite order by Lemma \ref{lemfo2}, and  
since $f'$ has finitely many zeros, applying a theorem of Bergweiler and Eremenko \cite{BE}
%together with its counterpart for finite lower order \cite{hinchliffe}, 
shows that $f$ has finitely many 
asymptotic values, each corresponding to a direct transcendental singularity of the inverse function \cite{BE,NEV}.

\begin{lem}
 \label{lemnot0}
Neither $0$ nor $\infty$ is an asymptotic value of $f$, and $f$ has infinitely many zeros and poles. 
\end{lem}
\begin{proof}
Let $g$ be $f$ or $1/f$ and assume that $\infty$ is an asymptotic value of $g$. Let $\sigma$ be a small positive real number and let $N_1$ and the set $E_1$ be as in Lemma \ref{lemfo4}. Since the inverse function $g^{-1}$ has finitely many singular values, Lemmas \ref{BELlem}
and \ref{lemfo4} give a large positive real number $N_2$ such that, for all large $r \not \in E_1$,
\begin{equation}
 \label{fo5}
|g(re^{i \theta } )| \leq N_2  \quad \hbox{for $\sigma \leq \pm \theta \leq \pi - \sigma$.}
\end{equation}

Since $\infty$ is by assumption an asymptotic value of $g$ and therefore a direct transcendental singularity
of $g^{-1}$ \cite{BE} there exists an unbounded component $U$ of the set $\{ z \in \C : |g(z)| > N_2 \}$ on which
$g$ has no poles, and the function
$$
u(z) = \log \left| \frac{g(z)}{N_2} \right| \quad (z \in U),  \quad
u(z) = 0 \quad (z \in \C \setminus U),
$$
is continuous, subharmonic and non-constant in the plane. Let $\theta^*(r)$ be defined
as in Lemma \ref{lemhm}. Then (\ref{fo5}) implies that $\theta^*(r) \leq 4 \sigma$ for all large
$r \not \in E_1$, and $E_1$ has upper logarithmic density at most $1/2$. Let $r$ be large and let
$R \geq 4r$. Then (\ref{p3}) gives
$$
\log B(R, u) \geq \pi \int_{2r}^{R/2} \frac{ds}{s \theta^*(s)} - O(1)
\geq \frac{\pi}{4 \sigma} \int_{ [2r, R/2] \setminus E_1} \frac{ds}{s} - O(1) \geq 
\frac{\pi}{12 \sigma} \log R 
$$
%\begin{eqnarray*}
% \log B(R, u) &\geq& \pi \int_{2r}^{R/2} \frac{ds}{s \theta^*(s)} - O(1) \\
%&\geq& \frac{\pi}{4 \sigma} \int_{ [2r, R/2] \setminus E_1} \frac{ds}{s} - O(1) \geq 
%\frac{\pi}{12 \sigma} \log R 
%\end{eqnarray*}
as $R \to \infty$. Since $\sigma$ may be chosen arbitrarily small, it follows that $u$ has infinite order.
But Poisson's formula leads to
$$
B(R, u) \leq \frac3{2\pi} \int_0^{2\pi} u(2R e^{i\theta} ) \, d \theta 
 \leq 3 m(2R, f) + 3 m(2R, 1/f) + O(1),
$$
and $f$ has finite order, which gives an immediate contradiction. It remains only to observe that $f$ must have infinitely many zeros and poles, by Iversen's theorem.
\end{proof}

\begin{lem}
 \label{lemav}
There exists $\alpha \in \C \setminus \R$ such that the set of asymptotic values of $f$ is precisely $\{ \alpha, \overline{\alpha}\}$. Moreover $f$ cannot tend to both $\alpha$ and $\overline{\alpha}$ on paths tending to infinity in the same component of $\C \setminus \R$.
\end{lem}
\begin{proof} By Lemma \ref{lemnot0}, neither $0$ nor $\infty$ is an asymptotic value of $f$. Suppose that $f$ has exactly one asymptotic value $a$. Since $f$ has infinitely many poles, $h = (f-a)/f'$ is transcendental with finitely many poles, and by a result of Lewis, Rossi and Weitsman \cite{LRW} there exists a path $\gamma$ tending to infinity with
$$
\int_\gamma \left| \frac{f'(z)}{f(z)-a} \right| \, |dz| = 
\int_\gamma \left| \frac{1}{h(z)} \right| \, |dz| < \infty.
$$
But then integration shows that $f(z)-a$ tends to a non-zero finite value as $z$ tends to infinity on $\gamma$, a contradiction.

Now let $\varepsilon $ be small and positive, so small that $|a - a'| \geq 2 \varepsilon $ whenever $a$ and $a'$ are distinct singular values of $f^{-1}$. Let $a$ be an asymptotic value of $f$. Then there exists a component $C_a$ of the set $\{ z \in \C : |f(z)-a| < \varepsilon \}$ containing a path tending to infinity on which $f(z)$ tends to $a$. By a standard argument \cite[p.287]{NEV} the function 
$\phi(t) = f^{-1} (a + e^{-t } ) $
maps the half-plane $\hbox{Re} \, (t) > \log 1/\varepsilon $ univalently onto $C_a$. Furthermore, the component $C_a$ contains infinitely many paths $\gamma_{a, j} $, each tending to infinity and mapped by $f$ onto %the open line segment 
$L_a = \{ a + t: 0 < t < \varepsilon /2 \}$. If $a$ is non-real then the 
$\gamma_{a, j}$ do not meet $\R$, since $\varepsilon $ is small and $f$ is real. On the other hand, if $a$ is real and $\gamma_{a, j} $ meets $\R$ then $\gamma_{a, j} \subseteq \R$ since $f$ is real and has no critical values $w \in L_a$. Moreover if $a$ is real then $f(z)$ also tends to $a$ as $z$ tends to infinity on the path $\overline{\gamma_{a, j}}$.

It follows that if $f$ has at least three distinct asymptotic values, or at least two distinct real asymptotic values, then there exist disjoint simple paths $\lambda_1$ and $\lambda_2$ tending to infinity with the following properties. Either both paths 
%start at some $x \in \R$ and, apart from their common starting point either both paths 
lie in the upper half-plane $H^+$, or both in the lower half-plane $H^-$, and as $z$ tends to infinity on $\lambda_j$ the function $f(z)$ tends to $b_j \in \C \setminus \{ 0 \}$ with $b_1 \neq b_2$. Choose a large positive real number $R_1$, in particular so large that $f$ has no non-real zeros $z$ with $|z| \geq R_1$. 
This gives an unbounded domain $D_1$ with no zeros of $f$ in its closure, bounded by a subpath of $\lambda_1$, a subpath of $\lambda_2$ and an arc of the circle $S(0, R_1)$. Since $b_1 \neq b_2$ the Phragm\'en-Lindel\"of principle forces $1/f$ to be unbounded on $D_1$, which implies the existence of a direct transcendental singularity of $f^{-1}$ over $0$, contradicting Lemma~\ref{lemnot0}.

Thus $f$ has exactly two distinct asymptotic values, of which at most one is real. Since $f$ is real the the set of asymptotic values of $f$ is $\{ \alpha, \overline{\alpha}\}$ for some non-real $\alpha$, and the last assertion of the lemma follows from the argument of the previous paragraph.
\end{proof}

\section{The multiplicities of the poles of $f$}\label{fpoles}

In this section it will be shown that all but finitely many poles of $f$ are simple.
It follows from Lemma \ref{lemav} and the fact that $f'$ has finitely many zeros that a simple closed polygonal path $J$ may be chosen with the following properties. $J$ is symmetric with respect to the real axis, and all non-real finite singular values of $f^{-1}$ lie on $J$. Moreover, 
$J \cap \R = \{ -R, R \}$, where the positive real number $R$ is chosen so that all the finitely many real critical values of $f$ lie in the interval $[ - R/2, R/2 ]$, which in turns lies in the interior domain $D_1$ of $J$. Let $D_2$
be the complement of $J \cup D_1$ in $\C \cup \{ \infty \}$. Thus the $D_j$ are simply connected domains, with $0 \in D_1$ and $\infty \in D_2$.

Then each component $C$ of $f^{-1}(D_2)$ is simply connected, and contains exactly one pole of $f$ of multiplicity $p$, say, and is mapped $p:1$ onto $D_2$ by $f$. This follows from the fact that $f^{-1}$ has no singular values in $D_2 \setminus \{ \infty \}$,  and may be proved (see \cite[p.362]{La14} or \cite{Lanew}) by choosing a quasiconformal mapping $\phi$ which satisfies $\phi (D_2) = \{ w \in \C \cup \{ \infty \} : |w| > 1 \}$ and $\phi ( \infty ) = \infty$, and writing $\phi \circ f = g \circ \psi$, where $g$ is meromorphic and $\psi$ is quasiconformal, following which the argument from \cite[p.287]{NEV} is applied to~$g$. 

Next, let $D_3 = D_1 \setminus (-R, R/2]$. Then $D_3$ is a simply connected domain containing no singular values of $f^{-1}$, and all components of $f^{-1}(D_3)$ are conformally equivalent to $D_3$ under $f$. Since $f'$ has finitely many zeros a standard argument \cite[Lemma 4.2]{BEL} then shows that each component $C$ of $f^{-1} (D_1)$ contains finitely many components of $f^{-1}(D_3)$, so that $f$ is finite-valent on $C$. Moreover all but finitely many components of $f^{-1}(D_1)$ are conformally equivalent to $D_1$ under~$f$.
These considerations lead at once to the following lemma.

\begin{lem}\label{lemcomp}
All but finitely many poles $w$ of $f$ lie in
components $C= C_w$ of $f^{-1}(D_2)$ such that $C$ and all components $D$ of $f^{-1} (D_1)$ with $\partial C \cap \partial D \neq \emptyset $ have the following properties:\\
(a) $f'$ has no zeros in $C \cup \partial C$, and $w$ is the only pole of $f$ in $C$;\\
(b) $f'$ has no zeros in $D \cup \partial D$;\\
(c) the mapping $f: D \to D_1$ is a conformal bijection;\\
(d) $f$ has no non-real zeros in $D$.
\end{lem}

The next lemma is the key to the proof of Theorem \ref{thm1}.

\begin{lem}
\label{lemsimple}
All but finitely many poles of $f$ are simple. 
\end{lem}
\begin{proof} Assume that $f$ has infinitely many multiple poles. Then $f$ has a multiple pole $w$ satisfying the conclusions of Lemma \ref{lemcomp}. Let $C = C_w$ be the component of $f^{-1}(D_2)$ which contains $w$. Then $C$ is simply connected and $\partial C$ consists of finitely many pairwise disjoint piecewise smooth simple curves, each tending to infinity in both directions. Let $\Gamma$ be a component of $\partial C$. Then $\Gamma \subseteq \partial D$ for some component $D$ of $f^{-1}(D_1)$, and $f$ is univalent on $D$ and so on $\Gamma$ by Lemma \ref{lemcomp}(c). As $z$ tends to infinity along $\Gamma$ in each direction, the image $f(z)$ tends to either $\alpha$ or $\overline{\alpha}$. Here $\Gamma$ will be called a type A component of $\partial C$ if $f(\Gamma)$ is a component of $J \setminus \{ \alpha, \overline{\alpha} \}$. If this is not the case then $f(\Gamma)$ is either $J \setminus \{ \alpha \}$ or $J \setminus \{ \overline{\alpha} \}$, and $\Gamma$ will be called type~B.

%The reflection of a component $\Gamma$ of $\partial C$ across the real axis is also a component of %$\partial C$, since $f$ is real, and either $\Gamma = \overline{\Gamma}$ or $\Gamma \cap \overline{\Gamma} %= \emptyset$. 
Every type A component $\Gamma$ of $\partial C$ must meet the real axis, by Lemma \ref{lemav}, and so is symmetric with respect to $\R$, since $f$ is real. Conversely, a component $\Gamma$ of $\partial C$ which meets  the real
axis has $\Gamma = \overline{\Gamma}$ and must be type A since the real function $f$ then has asymptotic values $\alpha$ and $\overline{\alpha}$ on $\Gamma$. 

If $\Gamma$ is a type B component of $\partial C$ and $D$ is that component of $f^{-1}(D_1)$ which satisfies
$\partial D \cap \Gamma \neq \emptyset$ then $\Gamma = \partial C$ by Lemma \ref{lemcomp}(b) and (c).
Thus at least one component of $\partial C$ must be type A, and hence at least two. On the other hand $\partial C$ cannot have three type A components, since one would have to separate the other two. It follows that since $w$ is a multiple pole there must be exactly two type A components $A_1$ and $A_2$ of $\partial C$, and at least one type B component $B_1$. Here $B_1$ is the boundary of a component $E_1$ of $f^{-1} (D_1)$. This component $E_1$ cannot meet the real axis since otherwise the fact that $B_1 \cap \R = \emptyset$ gives $\R \subseteq E_1$, whereas $A_q \cap \R \neq \emptyset$ for $q = 1, 2$. Now $E_1$ must contain a zero of $f$, but cannot contain a non-real zero of $f$ by Lemma \ref{lemcomp}(d). This is a contradiction and Lemma \ref{lemsimple} is proved.
\end{proof}

\section{Completion of the proof of Theorem \ref{thm1}}

It now follows from Lemma \ref{lemfo2}, Lemma \ref{lemsimple} and the fact that $f'$ has finitely many zeros that the Schwarzian derivative
\begin{equation}
S_f = \frac{f'''}{f'} - \frac{3}{2} \left( \frac{f''}{f'} \right)^2 = 2B
\label{s0}
\end{equation}
of $f$ is a rational function, and is not identically zero since $f$ is transcendental. Let
\begin{equation}
 \label{s1}
B(z) = b^2 z^n (1+ o(1)) \quad \hbox{as $z \to \infty$, where $b \in \C \setminus \{ 0 \}$ and $n \in \Z$.}
\end{equation}
Since $f^{-1}$ has at least two direct transcendental singularities the order of $f$ is at least $1$, and so (\ref{s0}) and an application of the Wiman-Valiron theory \cite{Hay5} to $1/f'$ show that $n \geq 0$ in (\ref{s1}). 

The following argument is self-contained but uses some methods similar to those of \cite{nicks}.
If $u$ and $v$ are linearly independent solutions of the differential equation
\begin{equation}
 \label{s2}
w'' + B(z) w = 0
\end{equation}
on a simply connected domain $D$ on which $B$ has no poles, then there exists a M\"obius transformation $T_0$
such that $f = T_0(u/v)$ on $D$ \cite[Chapter 6]{Lai1}. Hence $u/v$ and $(u/v)'$ extend to be meromorphic on $\C$ and, since the Wronskian of $u$ and $v$ is constant, so do $v^2$, $u^2$ and $uv$.

The equation (\ref{s2}) has $n+2$ distinct critical rays \cite{Hil2}, namely those rays $\arg z = \theta $ such that
$$
2 \arg b + (n+2) \theta = 0 \quad \hbox{(mod $2 \pi$).}
$$
If $\arg z = \theta_0$ is a critical ray and $\varepsilon $ and $1/R_0$ are small and positive then in the sectorial region
$$
S(R_0, \theta_0, \varepsilon ) = \left\{ z \in \C : |z| > R_0, \quad | \arg z - \theta_0 | < \frac{2\pi}{n+2} - \varepsilon \right\} 
$$
there are principal solutions $u_1$, $u_2$ of (\ref{s2}) given by \cite{Hil2}
$$
u_j(z) \sim B(z)^{-1/4} \exp( (-1)^j i Z ), \quad
Z = \int_{R_0 e^{i\theta_0}}^z B(t)^{1/2} \, dt \sim \left( \frac{2 b}{n+2} \right) z^{(n+2)/2} .
$$

\begin{lem}
 \label{lemcr}
The critical rays of (\ref{s2}) are the positive and negative real axis.
\end{lem}
\begin{proof} Suppose that $\arg z = \theta_0$ is a critical ray, where $\theta_0 \in (0, \pi) \cup (\pi, 2 \pi)$. Let $ \delta $ be small and positive. Then the $u_j$ are such that, without loss of generality, 
$u_1(z)/u_2(z)$ tends to infinity on the ray $\arg z = \theta_0 + \delta$ and to zero on the ray 
$\arg z = \theta_0 - \delta$. Since $f = T_1(u_1/u_2)$ for some M\"obius transformation $T_1$, this gives distinct asymptotic values of $f$ approached on paths in the same component of $\C \setminus \R$, contradicting Lemma \ref{lemav}. 
\end{proof}

It now follows that $n = 0$ and that $b$ may be chosen to be real and positive in (\ref{s1}). 
%By a change of variables it may be assumed that $b = 1$ and that 
Thus (\ref{s2}) has principal solutions satisfying
$$
u_j(z) = \exp( (-1)^j i ( bz + O( \log |z| )) ) \quad  \hbox{in} \quad
\{ z \in \C : |z| > R_0, \quad | \arg z  | <  \pi - \varepsilon \} .
$$
Each function $w = w_j = u_j^2$ is meromorphic in the plane and by (\ref{s2}) satisfies 
$$
2 ww'' - (w')^2 + 4 B w^2 = 0 ,
$$
so that $w_j$ has finitely many poles in $\C$. Applying the Phragm\'en-Lindel\"of principle to the functions
$v_j(z) = u_j(z)^2 \exp( (-1)^{j+1} 2 i b z) $, which are of finite order in the plane with finitely many poles, then shows that $v_1$ and $v_2$ are rational functions. Hence
$$\frac{u_2(z)^2}{u_1(z)^2} = S_1(z) \exp( 4ibz) $$
with $S_1$ a rational function, and since $u_2/u_1$ is meromorphic in the plane it follows that there exists a rational function $R$ 
%and a M\"obius transformation $T$ 
such that $u_2(z)/u_1(z) = R(z) e^{2ibz}$.
%$$
%f = T\left(\frac{u_2}{u_1}\right), \quad \frac{u_2(z)}{u_1(z)} = R(z) e^{2ibz} .
%$$
Since $f$ has infinitely many zeros and poles by Lemma \ref{lemnot0}, it may be assumed that
$$
f(z) =  \frac{ R(z)e^{icz} - 1}{A R(z)e^{icz} + A'},
$$
where $c = 2b > 0$ and $A, A' \in \C \setminus \{ 0 \}$. But $f$ has infinitely many real zeros
and so $R(z) \overline{R( \overline{z})} \equiv 1$. Further, $f$ has asymptotic values $1/A$ and $-1/A'$
and Lemma \ref{lemav} gives $A' = - \overline{A} \not \in \R$, which proves (\ref{conc1}) and (\ref{conc}). 
The remaining assertions of Theorem \ref{thm1} now follow from Lemma \ref{lemelem2}.

\end{document}